\documentclass[12pt,amstex]{article}
\usepackage{amsmath,amssymb,latexsym}
\makeatletter

\@addtoreset{equation}{section}
\makeatother
\def\tag#1{\eqno{\rm{(#1)}}}

\begin{document}
\title{Analytic General Solutions of Nonlinear Difference Equations }
\author{Mami SUZUKI\thanks{Research partially supported by the Grant--in--Aid 
for Scientific Research (C) 15540217 from the Ministry of Education, 
Science and Culture of Japan.}\\
}
\date{
Department of Management Informatics, \\
Aichi Gakusen Univ. Japan.\\
e-mail: m-suzuki@gakusen.ac.jp, 
}
\maketitle

\begin{abstract}
There is no general existence theorem for solutions for nonlinear difference equations, 
so we must prove the existence of solutions in accordance with models one by one.

In our work, we found theorems for the existence of 
analytic solutions of nonlinear second order difference equations. The main work of 
the present paper is  
obtaining representations of analytic general solutions with new 
methods of complex analysis. 
\end{abstract}
\quad\\
Subj-class: \,\,Classical Analysis and ODEs.\\
MSC-class:\,\,39A10,39A11,39B32.\par

\section{Introduction}

\vspace{0.7cm}
There is a general existence theorem for solutions of analytic differential 
equations, but 
we have no general existence theorem for 
analytic difference equations. 
For example, we consider the following first order nonlinear difference equation 
$$  x(t+1) = 2 x(t) + x(t)^2. \eqno{(*)}  $$ 
Putting $x(t) = -1 + y(t),$ we get $y(t+1) = y(t)^2$ and $\log y(t+1) = 2 \log y(t).$ 
Then $u(t) = 2^t$ is a (particular) solution of the equation $u(t+1) = 2u(t).$ 
Putting $C(t) = (\log y(t))/u(t),$ we have $C(t+1) = C(t),$ that is $C(t) = \pi(t),$ where 
$\pi(t)$ is an entire solution with the period $1.$ 
Therefore, a general solution of $(*),$ which tends to $0$ as 
$t \to - \infty,$ is given by $x(t) = \exp[\pi(t) 2^t] - 1.$ 

This simple example $(*)$ illustrates the whole make-up of the present paper. 
$0$ is its equilibrium point of $(*)$, with the characteristic value $2.$ 
A formal solution of it is obtained by putting $x(t) = \sum_{n=1}^{\infty} a_n (2^t)^n.$  
If its convergence is shown, then we have a solution $x(t)$ of the initial value problem, 
which tends to $0$ as $t \to - \infty.$ Further we proceed to seek general solutions. 

For analytic differential equations, a solution of initial value problem is always represented 
by a power series. 
This is the reason that the general existence theorem can be established for differential equations. 
But for difference equations, this is not the case. Next we consider the following 
first order difference equation  
$$  x(t+1) = x(t) + x(t)^2,  \eqno{(**)} $$ 
for which $0$ is the equilibrium point with characteristic value $1,$ 
but we can not put its formal solution in the form 
$\sum_{n=1}^{\infty} a_n (1^t)^n.$ That is, selection of appropriate formal
 solution depends on the problem. 
Of course, by \cite{Kimu} p.237 Theorem 14.2, 
($**$) has a local solution with the asymptotic expansion
$$ x(t) \sim -\frac{1}{t} \left\{1 + \sum_{j+k \geq 1} \hat{q}_{jk} t^{-j} 
\left(\frac{\log t}{t} \right)^k \right\}^{-1}, $$ 
where $\hat{q}_{jk}$ are constants.
\par
 
Further, for analytic differential equations, the solution is determined uniquely 
by the initial condition. However, 
for analytic difference equations, solution cannot be determined by 
the (initial) condition $x(t) \to 0$ as $t \to - \infty,$ 
hence we need to consider general solutions.

In this paper, we consider the following second order nonlinear difference equation,
\begin{equation}
u(t+2)=f(u(t),u(t+1)), \label{1.1}
\end{equation}
where $f(x,y)$ is an entire function of $x$, $y$.  
We assume that there is an equilibrium point $u^*: u^* = f(u^*, u^*).$ 
We can take $u^* = 0,$ that is $f(0,0)= 0$ without losing generality. \par

Many studies for difference equations are considered with discrete variables. 
Indeed such the equation (\ref{1.1}) is often considered for $t\in {\Bbb N}$. 
However in our study, 
we consider the difference equation (\ref{1.1}) with a continuous variable $t$. 
If "$t$" of equation (\ref{1.1}) represents "time", then $t$ is of course 
a real variable. However hereafter, in (\ref{1.1}), $t$ represents a complex variable, 
because we consist more general theorems. \par 
 
Our aim is to obtain analytic general solutions $u(t)$ of (\ref{1.1}) such that 
$u(t+n) \to 0$ as 
$n \to +\infty$ or $n \to - \infty.$  \par

We define $f(x,y)$ in (\ref{1.1}) such that   
\begin{equation}
 f(x,y) = - \beta x - \alpha y + g(x,y), \quad \beta \ne 0,  \label{1.2}
\end{equation}
where $g$ consists of higher order terms for $x$, $y$ such that 
$g(x,y)=\sum_{i,j\geqq 0, i+j\geqq 2} b_{i,j} x^iy^j\not\equiv 0,$ 
and $\alpha$, $\beta$, $b_{i,j}$ are constants.    
Further we assume that at least one of moduli of the  
characteristic values is neither $0$ nor $1.$ 
The case that both of characteristics equal to $1$ will be treated in another 
paper. \par
\quad

The processes of my work are as follows: 
{\bf 1)} determination of formal solutions, {\bf 2)} getting particular solution by Schauder's 
Fixed Point Theorem in a locally convex topological space, 
{\bf 3)} obtaining general solutions by the method of Kimura\cite{Kimu}
and Yanagihara\cite{Ya}.  \par


\section{Analytic Solutions}


\subsection{ A formal solution.}

 The characteristic equation of (\ref{1.1}) with (\ref{1.2}) is   
\begin{equation}
D(\lambda)=\lambda^2+\alpha\lambda+\beta=0.\label{2.1}
\end{equation}
 Let $\lambda_1$, $\lambda_2$ be roots of the characteristic equation and 
$|\lambda_1|\leqq |\lambda_2|$. 
Then we consider following two cases, i) $|\lambda_1|< 1$ and ii) $|\lambda_2|>1$. 
Of course, some characteristic equations have properties both 
i) and ii). \par

In case i), we consider solutions such that 
$$u(t+n)\to \,0,\quad \text {as}\,\,\,n\to\, +\infty.$$ 
In case ii), we consider solutions such that 
$$u(t+n)\to \,0,\quad \text {as}\,\,\,n\to\, -\infty.$$

In case i) we put $\lambda=\lambda_1$, and in case ii) we put 
$\lambda=\lambda_2$. Then we
put a formal solution to (\ref{1.1})
$$u(t)=\sum_{n=1}^{\infty} a_n\lambda^{nt},$$
in both cases.
We substitute $u(t)=\sum_{n=1}^{\infty} a_n\lambda^{nt}$, 
$u(t+1)=\sum_{n=1}^{\infty} a_n\lambda^{n(t+1)}$, 
$u(t+2)=\sum_{n=1}^{\infty} a_n\lambda^{n(t+2)}$, 
into (1.1). And we compare the coefficients of $\lambda^{nt},\,(n=1,2,\cdots)$, 
then we have, with $D(\lambda)$ in (\ref{2.1}),
\begin{equation}
\left\{
\begin{array}{l}
a_1\cdot D(\lambda)=0,\notag\\
a_2\cdot D(\lambda^2)=
a_1^2(b_{2,0}+b_{1,1} \lambda+ b_{0,2}\lambda^2),
\\
a_3\cdot D(\lambda^3)=
b_{2,0}2 a_1a_2+b_{1,1} a_1 a_2\lambda(\lambda+1)
+b_{0,2}2 a_1 a_2\lambda^3\notag\\
\qquad\qquad\qquad\qquad\qquad\qquad\qquad\qquad\quad\,\,\,\,\,+
a_1^3(b_{3,0}+b_{2,1}\lambda+b_{1,2}\lambda^2+b_{0,3}\lambda^3),\notag
\\
\cdots\\
a_k\cdot D(\lambda^k)=C_k(a_1,\cdots,a_{k-1}),\quad\notag\\
\cdots,\notag
\end{array}
\right.%
\end{equation}
where $C_k(a_1,\cdots,a_{k-1})$ are 
polynomials of $a_1,\cdots,a_{k-1}$ with coefficients 
$b_{i,j}\lambda^l$,  
$0\leqq i \leqq k$, $0\leqq j\leqq k$, $0\leqq l\leqq k$, $2\leqq i+j\leqq 
k$.
From definition of $\lambda$ and $D$, we have 
$D(\lambda)=0$ and $D(\lambda^k)\neq 0\, (k\geqq 2)$, and we can have $a_1$ is arbitrary. \par 
Here we suppose that $a_1\neq 0$. Then we have 
\begin{equation}
a_k=\frac{a_1^k}{D(\lambda^k)}C^*_k(b_{i,j},\lambda^l), \,\, k\geqq 2, \label{2.2}
\end{equation}
where $C^*_k(b_{i,j},\lambda^l)$ are constants which are given by the function $f$, in which
they consist of $b_{i,j}$, $2\leqq i+j\leqq k$ and $\lambda^l$, $0\leqq l\leqq k$. 
Hence we can determine a formal solution of (\ref{1.1}), 
\begin{equation}
u(t)=\sum_{n=1}^{\infty} a_n\lambda^{nt},\label{2.3}
\end{equation} 
in both cases i) and ii). Here we have $a_1$ is arbitrary not $0$,  
and $a_k$ are determined by $a_1$. 


\subsection{ Existence of an analytic solution }
Here we put $u(t)=s,u(t+1)=w, u(t+2)=z$, and 
$H(s,w,z)=-z+f(s,w)$. 
%
Then the equation (\ref{1.1}) can be written such as
\begin{equation}
H(u(t),u(t+1),u(t+2))=0.\label{2.4}
\end{equation}
$H(s,w,z)$ is holomorphic in a neighborhood of $(0,0,0)$ and 
we have $H(0,0,0)=0$ easily. 
Furthermore we have 
$
\frac{\partial H}{\partial s}(0,0,0)
=\frac{\partial f}{\partial s} \Bigr |_{s=w=0}
=-\beta \neq 0$ as remarked in (\ref{1.2}). 
From implicit function theorem, 
 for the equation $H(s,w,z)=0$, 
we have a holomorphic function 
$\phi$ such that 
\begin{equation}
s=\phi(w,z) \quad \mbox{for}\quad |w|,\,|z|\leqq \rho\label {2.5}
\end{equation}
for some $\rho>0$. 
Furthermore we have a constant $K$ such that 
\begin{equation}
|s|=|\phi(w,z)|\leqq K(|w|+|z|) \quad \mbox{for}\quad |w|,|z|\leqq \rho.\label{2.6}
\end{equation}
\quad
Let $N$ be a positive integer. Put the partial sum of formal solution
as 
$P_N(t)=\sum_{n=1}^N \alpha_n \lambda^{nt}$, and put 
$p_N(t)=u(t)-P_N(t)$. Here we rewrite $p(t)=p_N(t)$. \par
Moreover we define following sets,
\begin{align}
&S(\eta)=\{t\in \Bbb{C}: |\lambda^t|\leqq \eta \},\notag\\
&J(A,\eta)=\{p:p(t) \,\text{is holomorphic and }\, |p(t)|\leqq A|\lambda^t|^{N+1} \,
\mbox{for }\, t\in S(\eta)\}\notag.
\end{align}
in which $A>0$ and $\eta$, $0<\eta<1$ are constants. We determined these constants in a proof 
of existence for a fixed point of following maps $T_i$ ($i=1,2$).\par
 Suppose there would exist a solution $u(t)$ of (\ref{1.1}) in $S(\eta)$. 
Then $p_N(t)=u(t)-P_N(t)$ would 
belong to $J(A,\eta)$ for some suitably chosen constants $A$, $\eta$, and would satisfy the equation
\begin{equation}
p(t+2)=f(p(t)+P_N(t),p(t+1)+P_N(t+1))-P_N(t+2),\label{2.7}
\end{equation}
with $p(t)=p_N(t)$. 
Conversely, suppose there would exist a solution $p(t)$ of (\ref{2.7}), then $u(t)=p(t)+P_N(t)$ 
would be a solution of (\ref{1.1}). 
So, hereafter we concentrate on 
proving the existence of $p(t)\in J(A,\eta)$ such that $u(t)=p(t)+P_N(t)$ satisfies (\ref{2.7}).

In case i) $|\lambda|<1$, 
the existence of solutions $u(t)$ of (\ref{2.7}) is equivalent to 
the existence of $p(t)$ which satisfies    
$$p(t)
=\phi(p(t+1)+P_N(t+1),p(t+2)+P_N(t+2))-P_N(t).$$
For $p(t)\in J(A,\eta)$, we put 
\begin{equation}
T_1[p](t)=\phi(p(t+1)+P_N(t+1),p(t+2)+P_N(t+2))-P_N(t).\label{2.8}
\end{equation}
Then we can prove that $T_1$ maps $J(A,\eta)$ into itself (see Appendix A). 
The map $T_1$ is obviously continuous if $J(A,\eta)$ is 
endowed with topology of uniform 
convergence on compact sets in $S(\eta)$. Furthermore 
$J(A,\eta)$ is clearly convex, and is relatively compact by the theorem of 
Montel \cite{Ahlfors}.\par

By Schauder's fixed point theorem \cite{Dug}(p.74), \cite{Smar}(p.32), we obtain the existence of a 
fixed point $p(t)=p_N(t)\in J(A,\eta)$ of $T_1$ in $S(\eta)$. 
Moreover we can prove uniqueness of the fixed point (see Appendix B) and independence from $N$ 
(see Appendix C). 
Hence 
we have an analytic solution $u(t)$ in $S(\eta)$. 

\par
\vspace{0.5cm}


\quad
In case ii) $|\lambda|>1$, 
(\ref{2.7}) is equivalent to 
the existence of $p(t)$ which satisfies,   
$$p(t)=f(p(t-2)+P_N(t-2),p(t-1)+P_N(t-1))-P_N(t).$$ 
For $p(t)\in J(A,\eta)$, we put
\begin{equation}
T_2[p](t)=f(p(t-2)+P_N(t-2),p(t-1)+P_N(t-1))-P_N(t)\notag
\end{equation}
Then we can prove the existence of an analytic solution $u(t)$ in 
$S(\eta)$  by the arguments similar as above. 
\par
\vspace{0.7cm}


Thus we have the following Theorem 1.
\par
\vspace{0.7cm}


 {\bf Theorem 1.} \em
Let $\lambda_1,\, \lambda_2$ be roots of 
$D(\lambda)=0$ in (2.1), with $|\lambda_1|\leqq|\lambda_2|$.
Suppose  $|\lambda_1|<1$ or $|\lambda_2|>1$. 
Put $\lambda=\lambda_1$ for the former, and $\lambda=\lambda_2$ for latter. 
And we assume that $\lambda_1^k\neq \lambda_2$ and $\lambda_2^k\neq \lambda_1$ 
for any $k\in {\mathbb N}$. 
Then there is an $\eta>0$ such that we have a holomorphic 
solution $u(t)=\sum^{\infty}_{n=1}a_n\lambda^{nt}$ in 
$S(\eta)=\{t;|\lambda^t|<\eta\}$.
\em\par
\vspace{0.5cm}

{\bf In case ii).} The solution $u(t)$ can be analytically continued 
to the whole plane, by making use of the equation (\ref{1.1}) 
$u(t+2)=f(u(t),u(t+1)).$\par
{\bf In case i).} The function $\phi(w,z)$ in (\ref{2.5}) 
$s=\phi(w,z)$ for $|w|,\,|z|\leqq \rho$  
is defined only locally, 
though we can also analytically continue 
$u(t)$, keeping out of branch points. The solution obtained is 
multi-valued. \par
\vspace{0.5cm}

\subsection{Particular solutions.} 
In this subsection, we consider solutions $u_1(t)$ and $u_2(t)$ which are 
respectively depend on $\lambda_1$ and $\lambda_2$. Put formal solutions such that 
$u_1(t) = \sum_{n=1}^{\infty} a_{1,n} \lambda_1^{nt} $ and 
$u_2(t) = \sum_{n=1}^{\infty} a_{2,n} \lambda_2^{nt}$, we have  
$a_{m,k}\cdot D(\lambda_m^k)=C_{m,k}(a_{m,1},\cdots,a_{m,k-1}),$ 
$(m=1,2; k\in \Bbb{N})$ 
with the similar arguments in subsection {\bf 2.1}, where 
$C_{m,k}(a_{m,1},\cdots,a_{m,k-1})$ are 
polynomials of $a_{m,1},\cdots,a_{m,k-1}$ with coefficients 
$b_{i,j}\lambda_m^l$,  
$0\leqq i \leqq k$, $0\leqq j\leqq k$, $0\leqq l\leqq k$, $2\leqq i+j\leqq 
k$. 
Furthermore if we take as $a_{m,1}\neq 0$, then we have 
\begin{equation}
a_{m,k}D(\lambda_m^k)=a_{m,1}^k
C^*_{m,k}(b_{i,j},\lambda_m^l), \,\, m=1,2;\,\,k\geqq 2, \label{2.9}
\end{equation}
where $C^{*}_{m,k}(b_{i,j},\lambda_m^l)$ are constants which are given by the function $f$, 
in which 
they consist of $b_{i,j}$, $2\leqq i+j\leqq k$ and $\lambda_m^l$, $0\leqq l\leqq k$. 
Then we have the following lemma 2 and lemma 3 with the similar arguments 
in {\bf 2.1-2.2}.\par


\vspace{0.7cm}
{\bf Lemma 2.} \em 
 Let $\lambda_1, \lambda_2$ be roots of (\ref{2.1}) with $|\lambda_1| \leqq |\lambda_2| < 1$. 
If $\lambda_2^k \ne \lambda_1$ 
for any positive integer $k$ greater than $1$, 
then there are constants  
$\eta_1,\, \eta_2>0$ such that we have following 
two holomorphic solutions $u_1$ and $u_2$ of (\ref{1.1}), 
\begin{gather}
u_m(t) = \sum_{n=1}^{\infty} a_{m,n} \lambda_m^{nt} \quad \text{in} 
\quad S(\eta_m)=\{t;|\lambda_m^t|<\eta_m\}, \quad (m=1,2),\notag
\end{gather}
in which $a_{1,1}$ and $a_{2,1}$ can be taken to be arbitrary 
non-zero constants.\par 
For the case $\lambda_2^k = \lambda_1$ for some $k \in {\mathbb N},$ 
if $C^*_{2,k}(b_{i,j},\lambda_2^l)=0$ given in (\ref{2.9}), 
then we take $a_{2,1}\neq 0$ and $a_{2,k} \ne 0$ arbitrary, 
and have the solution $u_2(t)$ as above.   
On the other hand, if $C^*_{2,k}(b_{i,j},\lambda_2^l)\neq 0$ for the $k$, 
then we take $a_{2,j}=0$ for $(j\ne kn, n\in \Bbb{N})$, and can take   
$a_{2,k} \ne 0$ arbitrary, then we determine coefficients  
 $a_{2,kn}$ 
for $n \geqq 2$ as above. 
Hence then there is an 
$\eta_2>0$ such that we have holomorphic solutions $u_1$ and $u_2$ of (\ref{1.1}), 
$$u_2(t) = \sum_{n=1}^{\infty} a_{2,kn} \lambda_2^{knt} \quad\text{in} \,\,S(\eta_2),$$  
as well as $u_1(t) = \sum_{n=1}^{\infty} a_{1,n} \lambda_1^{nt}$ 
in $S(\eta_1). $ In the case of $\lambda_2^k=\lambda_1$ and 
$C^*_{2,k}(b_{i,j},\lambda_2^l)\neq 0$, 
if we take $a_{2,k}=a_{1,1}$, then $u_2(t)=u_1(t)$ in  $S(\eta_1)\cap S(\eta_2)$. \par
Thus, in the both cases, $u_1(t + n) \to 0, u_2(t + n) \to 0$ as $n \to \infty$ uniformly 
on any compact subset of the $t$-plane. 
\em\par
\vspace{0.5cm}

\quad {\bf Proof.}\quad 
If $\lambda_2^k \ne \lambda_1$ for any $k \in {\mathbb N},$ then we 
can determine formal 
solution 
$u_2(t) = \sum_{n=1}^{\infty} a_{2,n} \lambda_2^{nt}$ as in subsection {\bf 2.1}, 
with $\lambda = \lambda_2$ instead of $\lambda_1.$ And we can show that it is 
an actual solution  
as in subsections {\bf 2.1-2.2} \par
For the case $\lambda_2^k = \lambda_1$ for some $k \in {\mathbb N},$ 
if we take $a_{2,1}\neq 0$, form (\ref{2.9}), we have 
\begin{equation}
a_{2,k}D(\lambda_2^k)=a_{2,k}D(\lambda_1)=a_{2,1}^k
C^*_{2,k}(b_{i,j},\lambda_2^l)=0.\label{2.10}
\end{equation}
If $C^*_{2,k}(b_{i,j},\lambda_2^l)= 0$, we can take $a_{2,k}$ arbitrary, and 
determine $a_{2,n}$, $2\leqq n \leqq k-1, n\geqq k+1$ by $a_{2,1}$ as in (\ref{2.9}).\par
However if $C^*_{2,k}(b_{i,j},\lambda_2^l)\neq 0$, the equation (\ref{2.10}) 
is contradiction. Thus 
we must take $a_{2,1}=0$, then $a_{2,n}=0,$ for $n\leqq k-1$ by 
$a_{2,k}\cdot D(\lambda_2^k)=C_{2,k}(a_{2,1},\cdots,a_{2,k-1}).$ 
Then we can take $a_{2,k}$ to be 
arbitrary non-zero constant, and determine coefficients $a_{2,n}$ as follows 
\begin{equation}
a_{2,n}=
\begin{cases}
 0,\quad (n\neq km , \,\, m\in \Bbb{N}),\notag\\
a_{2,k}^m\frac{C_{2,m}^*(b_{i,j}, \lambda_2^{lk})}{D(\lambda_2^{km})}
=a_{2,k}^m\frac{C_{2,m}^*(b_{i,j}, \lambda_1^{l})}{D(\lambda_1^m)},
\quad (n= km , \,\, m\in \Bbb{N}),\notag
\end{cases}
\end{equation}
where $C^*_{2,m}(b_{i,j},\lambda_1^l)$ are constants defined in (\ref{2.9}). 
Hence we can determine a formal solution $u_2(t)$ such that 
$$u_2(t) = \sum_{n=1}^{\infty} a_{2,kn} \lambda_2^{knt} \quad\text{in} \,\,S(\eta_2).$$
If we take $a_{2,k}=a_{1,1}$, then we have only one solution.  
Futhermore for the both cases of $\lambda_2^k=\lambda_1$, 
we can prove that 
there is an $\eta_2>0$ such that we have a holomorphic solution 
$u_2=\sum_{n=1}^{\infty} a_{2,kn} \lambda_2^{knt}$ in $S(\eta_2)$, with the 
similar arguments in {\bf 2.2}. 
 in  $S(\eta_1)\cap S(\eta_2)$.
\par 
Obviously $u_1(t + n) \to 0, u_2(t + n) \to 0$ as $n \to +\infty$ uniformly on 
any compact subset of the $t$-plane. 
$\square$\par

%
%

\vspace{0.7cm}


{\bf Lemma 3.} \em
 Let $\lambda_1, \lambda_2$ be roots of (\ref{2.1}) with $1<|\lambda_1| \leqq |\lambda_2|$. 
If $\lambda_1^k \ne \lambda_2$ 
for any positive integer $k$ greater than $1$, 
then there are constants  
$\eta_1,\, \eta_2>0$ such that we have following 
two holomorphic solutions $u_1$ and $u_2$ of (\ref{1.1}), 
\begin{gather}
u_m(t) = \sum_{n=1}^{\infty} a_{m,n} \lambda_m^{nt} \quad \text{in} 
\quad S(\eta_m)=\{t;|\lambda_m^t|<\eta_m\}, \quad (m=1,2),\notag
\end{gather}
in which $a_{1,1}$ and $a_{2,1}$ can be taken to be arbitrary 
non-zero constants.\par 
For the case $\lambda_1^k = \lambda_2$ for some $k \in {\mathbb N},$ 
if $C^*_{1,k}(b_{i,j},\lambda_1^l)=0$ given in (\ref{2.9}), 
then we take $a_{1,1}\neq 0$ and $a_{1,k}\neq 0$ arbitrary, and 
have the solution $u_1(t)$ as above.   
On the other hand, if $C^*_{1,k}(b_{i,j},\lambda_1^l)\neq 0$ for the $k$, 
then we take $a_{1,j}=0$ for $(j\ne kn, n\in \Bbb{N})$, and can  
$a_{1,k} \ne 0$ arbitrary, then we determine coefficients  
 $a_{1,kn}$ 
for $n \geqq 2$ as above. 
Hence then there is an 
$\eta_1>0$ such that we have holomorphic solutions $u_1$ and $u_2$ of (\ref{1.1}), 
$$u_1(t) = \sum_{n=1}^{\infty} a_{1,kn} \lambda_1^{knt} \quad\text{in} \,\,S(\eta_1),$$  
as well as $u_2(t) = \sum_{n=1}^{\infty} a_{2,n} \lambda_2^{nt}$ 
in $S(\eta_2). $ In 
the case of $\lambda_1^k=\lambda_2$ and $C^*_{1,k}(b_{i,j},\lambda_1^l)\neq 0$,  
if we take $a_{1,k}=a_{2,1}$, then $u_1(t)=u_2(t)$ in  $S(\eta_1)\cap S(\eta_2)$. \par  
Thus, in the both cases, $u_1(t - n) \to 0, u_2(t - n) \to 0$ as $n \to \infty$ uniformly 
on any compact subset of the $t$-plane. \par
\em\par
\vspace{0.5cm}
\quad {\bf Proof.}\quad 
We can prove with arguments similar as Lemma 2. $\square$\par
\vspace{0.7cm}

The analytic solutions $u_1$ and $u_2$ obtained in Lemmas 2-3 are
"Particular Solutions" of (\ref{1.1}). 


\section{Analytic General Solutions}
Analytic general solutions of nonlinear difference equations have been 
investigated, for example, by Harris \cite{Harris}, \cite{Harris2}, 
and others, but we can not make use of their method. 
Here we follow the method of 
Kimura \cite{Kimu} and Yanagihara \cite{Ya}, 
where general solutions of the first order difference equations are studied. \par
In this section we consider the following case, 
$$
|\lambda_1|<1<|\lambda_2|.
$$
For other cases, we study 
general solutions of 
the difference equation (1.1) in other papers. \par

\vspace{0.4cm}
For a linear second order difference equation, general solutions are 
written by two particular solutions of it. But for a nonlinear second 
order difference equation, in this case, 
general solutions which converge to 
an equilibrium point of the equation are written by one of two 
particular solutions $u_1$ or $u_2$ of the difference equation.\par
\vspace{0.7cm}
 Let $u(t)$ be a solution of (\ref{1.1}), and $w(t) = u(t+1).$ Then (\ref{1.1}) 
can be written as a system of simultaneous equations 
\begin{equation}
 \begin{pmatrix} u(t+1) \\ w(t+1) \end{pmatrix} = \begin{pmatrix} 0 & 1 \\ 
 - \beta & - \alpha \end{pmatrix} \begin{pmatrix} u(t) \\ w(t) \end{pmatrix} 
+ \begin{pmatrix} 0 \\ g(u(t), w(t)) \end{pmatrix}  
\label{3.1} 
\end{equation} 
Let $\lambda_1, \lambda_2$ be roots of the equation (\ref{2.1}) and 
$P = \begin{pmatrix} 1 & 1 \\ \lambda_1 & \lambda_2 \end{pmatrix}.$ 
Put 
\begin{equation}
  \begin{pmatrix} u \\ w \end{pmatrix} = P \begin{pmatrix} x \\ y \end{pmatrix}. 
\label{3.2} 
\end{equation}
From $\lambda_1\neq\lambda_2$, we can transform the coefficient matrix of linear terms 
of (\ref{3.1}) into diagonal form, i.e.,   
(\ref{3.1}) is transformed to a following system with respect to $x, y:$ 

\begin{equation}
 \left\{  \begin{aligned} 
  x(t + 1) &= \lambda_1 x(t) + \sum_{i+j \geq 2} c_{ij} x(t)^i y(t)^j = X(x(t), y(t)), \\
  y(t + 1) &= \lambda_2 y(t) + \sum_{i+j \geq 2} d_{ij} x(t)^i y(t)^j = Y(x(t), y(t)). 
  \end{aligned}   \right.  
\label{3.3}  
\end{equation}
On the other hand, let  
$Q = \begin{pmatrix} 1 & 1 \\ \lambda_2 & \lambda_1 \end{pmatrix}.$ 
Put 
\begin{equation}
  \begin{pmatrix} u \\ w \end{pmatrix} = Q \begin{pmatrix} x \\ y \end{pmatrix}. 
\label{3.4} 
\end{equation}
Then (\ref{3.1}) is transformed to a system with respect to $x, y:$ 
\begin{equation}
 \left\{  \begin{aligned} 
  x(t + 1) &= \lambda_2 x(t) + \sum_{i+j \geq 2} c'_{ij} x(t)^i y(t)^j = X'(x(t), y(t)), \\
  y(t + 1) &= \lambda_1 y(t) + \sum_{i+j \geq 2} d'_{ij} x(t)^i y(t)^j = Y'(x(t), y(t)). 
  \end{aligned}   \right.  
\label{3.5}  
\end{equation}

Then we will show the following Theorem 4.

\vspace{0.7cm}

{\bf Theorem 4.} \em
Let $\lambda_1,\, \lambda_2$ be roots of the characteristic 
equation of (\ref{1.1}) such that $|\lambda_1|<1< |\lambda_2|$. 
Suppose that $u_1(t)$ and $u_2(t)$ are solutions of (1.1) which have the expansions 
$u_1(t)=\sum^{\infty}_{n=1} a_{1,n}\lambda_1^{nt}$ in $S(\eta_1)=\{t;|\lambda^t|<\eta_1 \}$,      
$u_2(t)=\sum^{\infty}_{n=1} a_{2,n}\lambda_2^{nt}$ in $S(\eta_2)=\{t;|\lambda^t|<\eta_2 \}$ 
with some constants $\eta_1,\eta_2>0$. 
Further suppose that $\Upsilon(t)$ is an analytic solution of (\ref{1.1}) 
such that either $\Upsilon(t+n)\to 0$ as $n\to +\infty$ or $n\to -\infty$,
 uniformly on any compact subsets of $t$-plane. 
If the solution $\Upsilon$ of (\ref{1.1}) satisfies 
$\Upsilon(t+(-1)^{m-1}n)\to 0$, $(m=1,2)$ as $n\to +\infty$, 
then 
there is a periodic entire function $\pi_m(t),(\pi_m(t+1)=\pi_m(t))$, such that 
\begin{align}
\Upsilon(t)&=\frac{1}{\lambda_{m+1}-\lambda_m}( \lambda_{m+1}\sum^{\infty}_{n=1}
a_{m,n}\lambda_m^{n(t+\pi_m(t))}-\sum^{\infty}_{n=1}
a_{m,n}\lambda_m^{n(t+\pi_m(t)+1)})\notag\\
&\qquad+\Psi_m\Biggr(  
\frac{1}{\lambda_{m+1}-\lambda_m}( \lambda_{m+1}\sum^{\infty}_{n=1}
a_{m,n}\lambda_m^{n(t+\pi_m(t))}-\sum^{\infty}_{n=1}
a_{m,n}\lambda_m^{n(t+\pi_m(t)+1)})
\Biggr),
\label{3.6}
\end{align}
in $S(\eta_m)$, with the convention $\lambda_3$ means $\lambda_1$. 
Further we have 
$\frac{\Upsilon(t+1+(-1)^{m-1}n)}{\Upsilon(t+(-1)^{m-1}n)}\to \lambda_m$, ($m=1,2$), 
as $n\to +\infty$.  

When $m=1$, 
$\Psi_1$ is a solution of 
\begin{equation}
\Psi(X(x,\Psi(x)))=Y(x,\Psi(x)),\label{3.7}
\end{equation}
and 
when $m=2$, 
$\Psi_2$ is a solution of 
\begin{equation}
\Psi(X'(x,\Psi(x)))=Y'(x,\Psi(x)),\label{3.8}
\end{equation}
in which 
$X$, $Y$ are defined in (\ref{3.3}), and $X'$, $Y'$ are defined in (\ref{3.5}).
\par
Conversely, a function $\Upsilon(t)$ which is represented as in (\ref{3.6}) 
in $S(\eta_m)$ for some $\eta_m>0$, where $\pi_m(t)$ is a periodic function with 
the period one,  
is a solution of (\ref{1.1}) such that $\Upsilon(t+(-1)^{m-1}n)\to 0$ 
and 
$\frac{\Upsilon(t+1+(-1)^{m-1}n)}{\Upsilon(t+(-1)^{m-1}n)}\to \lambda_m$ as 
$n \to +\infty$ with $m=1,2$.
\par
\em\par
\medskip
\vspace{0.7cm}

{\bf Proof.}\quad 
At first we prove the case $m=1$.\par 
Let $u(t)$ be the solution of (\ref{1.1}) in the argument of Section 2. And suppose 
$\Upsilon(t)$ be a solution of (\ref{1.1}) such that $\Upsilon(t+n)\to 0$ 
as $n\to +\infty$ 
uniformly on any compact subsets of $t$-plane.\par  

At first we will consider the meaning of the functional equation (\ref{3.7}). \par
Suppose (\ref{3.3}) admits a solution $(x(t), y(t))$. 
If 
$\frac{dx}{dt}\neq 0$, then we can write $t=\psi(x)$ with a function $\psi$ 
in a neighborhood of $x_0=x(t_0)$, and we can write 
\begin{equation}
y(t)=y(\psi(x))=\Psi(x),\label{3.9}
\end{equation}
as far as $\frac{dx}{dt}\neq 0$. 
Then the function $\Psi$ satisfies the functional equation (\ref{3.7}).\par
Conversely we assume that a function $\Psi$ is a solution of 
the functional equation (\ref{3.7}). If the first order difference equation 
\begin{equation}
x(t+1)=X(x(t),\Psi(x(t))),\label{3.10}
\end{equation}
has a solution $x(t)$, then we put $y(t)=\Psi(x(t))$ and have a 
solution $(x(t), y(t))$ of (\ref{3.3}). 
From \cite{Ya} we see that the first order difference equation 
(\ref{3.10}) has an analytic solution.\par
This relation is a point of our method.\par
\vspace{0.7cm}

Put $\omega(t)=\Upsilon(t+1)$ and 
\begin{equation}
\begin{pmatrix}\chi\\ \nu \end{pmatrix}
=P^{-1}\begin{pmatrix} \Upsilon\\\omega\end{pmatrix}.\label{3.11}
\end{equation}
Then we have 
$\chi(t)=\frac{1}{\lambda_2-\lambda_1}(\lambda_2\Upsilon(t)-\omega(t))$. 
Since $\Upsilon(t+n)\to 0$ and $\omega(t+n)\to 0$ as $n\to \infty$, we have 
$\chi(t+n) \to 0$ as $n\to \infty$.\par

Let $u(t)$ be a solution given in Section 2,  
$$u(t)=\sum^{\infty}_{n=1}a_{1,n} \lambda^{nt}\qquad\qquad (\lambda=\lambda_1).$$
Then we can write by (\ref{3.11}), 
since $\lambda_1 = \lambda$ and $u(t)$ is a function of $\lambda^t,$
\begin{equation}
 x(t)=\frac{1}{\lambda_2-\lambda}(\lambda_2 u(t)-u(t+1))
=\frac{1}{\lambda_2-\lambda}
\Biggr ( \sum_{n=1}^{\infty}(\lambda_2a_{1,n}-a_{1,n}\lambda^n)(\lambda^t)^n
\Biggr )
=\Tilde{U}(\lambda^t),\label{3.12}
\end{equation}
where $\zeta = {\tilde U}(\tau)$ is a function of 
$\tau = \lambda^t$ and ${\tilde U}'(0) = a_{1,1} \ne 0$ and ${\tilde U}(0) = 0.$ 
Since ${\tilde U}(\tau)$ is an open map, for any $\eta_1 > 0$ 
there is an $\eta_2 > 0$ such that 
$$ {\tilde U}(\{ |\tau| < \eta_1 \}) \supset \{|\zeta| < \eta_2 \}. $$
Since $\chi(t+n) \to 0$ as $n \to \infty,$ supposed that 
$t$ belongs to a compact set $K,$ there is an $n_0 \in {\mathbb N}$ 
such that for $t' \in K$  
$$|\chi(t'+n)|<\eta_2\quad (n\geqq n_0).$$ 
Thus there is a $\tau'=\lambda^{\sigma}$, such that
\begin{equation}
\chi(t'+n)=\Tilde{U}(\tau')=\Tilde{U}(\lambda^{\sigma}).\label{3.13}
\end{equation}
 Since $\Tilde{U}'(0)=a_{1,1}\neq 0$, using the theorem on implicit function 
we have a $\Tilde{U}^{-1}$ such that 
$$\lambda^{\sigma}=\Tilde{U}^{-1}(\chi(t'+n)).$$ 
Put $t=t'+n$, then $\lambda^{\sigma}=\Tilde{U}^{-1}(\chi(t))$, 
and we write
\begin{equation}
\sigma=\log_{\lambda}\Tilde{U}^{-1}(\chi(t))=\ell(t).\label{3.14}
\end{equation}
\quad When there is a solution $\chi(t)$ of (\ref{3.3}), 
from (\ref{3.10}), (\ref{3.12}-\ref{3.13}) 
we have 
\begin{align}
\chi(t+1)
&=X(\chi(t),\Psi(\chi(t)) )\notag\\
&=X(\Tilde{U}(\lambda^{\sigma}),\Psi(\Tilde{U}(\lambda^{\sigma}))  )\notag\\
&=X(x(\sigma),\Psi(x(\sigma)))\notag\\
&=x(\sigma+1)=\Tilde{U}(\lambda^{\sigma+1}).\notag  
\end{align}
Hence  
$$\sigma+1=\ell(t+1),\,\,\ell(t)+1=\ell(t+1).$$ 
If we put $\pi(t) = \ell(t) - t,$ then we obtain 
$\pi(t+1) = \ell(t+1) - (t+1) = \ell(t) - t = \pi(t),$ and we can write as  
\begin{equation}
 \ell(t) = t + \pi(t), \label{3.15}  
\end{equation}
where $\pi(t)$ is defined for a compact set $K$ with $\Re[t]$ 
sufficiently large. Furthermore we can continue the $\pi(t)$ analytically as a 
periodic function with the period $1.$ Thus we have  
$$ \sigma = t + \pi(t).  $$ 
From (\ref{3.13}) and (\ref{3.12}), $\chi(t)$ can be written as 
$$
\chi(t)=\Tilde{U}(\lambda^{t+\pi(t)})=x(t+\pi(t))
=\frac{1}{\lambda_2-\lambda_1}
(\lambda_2u(t+\pi(t))-u(t+1+\pi(t))).$$
We have following equations, making use of the equation (\ref{3.11}) 
\begin{align}
\Upsilon(t)&=\chi(t)+\nu(t)\notag\\
&=\chi(t)+\Psi(\chi(t))\notag\\
&=x(t+\pi(t))+\Psi(x(t+\pi(t)))\notag\\
&=\frac{1}{\lambda_2-\lambda_1}(\lambda_2\sum^{\infty}_{n=1}
a_{1,n}\lambda^{n(t+\pi(t))}-\sum^{\infty}_{n=1}
a_{1,n}\lambda^{n(t+\pi(t)+1)})\notag\\
&\qquad+\Psi\Biggr(  
\frac{1}{\lambda_2-\lambda_1}(\lambda_2 \sum^{\infty}_{n=1}
a_{1,n}\lambda^{n(t+\pi(t))}-\sum^{\infty}_{n=1}
a_{1,n}\lambda^{n(t+\pi(t)+1)})
\Biggr),\notag
\end{align}
where $\pi(t)$ is defined for $t\in\cup_{n\in \Bbb{Z}}(K + n)$ 
with a compact set $K.$ Since $K$ is arbitrary, 
we can continue $\pi(t)$ analytically to 
a periodic entire function with period $1,$ 
and 
$\Psi$ is a solution of (\ref{3.7}). By making use of the Theorem in  \cite{Suzu1}, (\cite{Suzu3}), 
$\Psi$ is obtained in the form, in a neighborhood of $x=0$, 
\begin{equation}
\Psi(x)=\sum_{n=2}^{\infty} \gamma_n x^n,\label{3.16}
\end{equation}
that is, the expansion begins with $x^2$. 
From $\chi(t+1)=X(\chi(t),\Psi(\chi(t)))$, we have 
$$\chi(t+1)=\lambda_1\chi(t)+\sum_{i+j\geqq 2}c_{ij}\chi(t)^i
\Psi(\chi(t))^j,$$
and
$$\frac{\chi(t+1)}{\chi(t)}=\lambda_1+\sum_{i+j\geqq 2}c_{ij}
\chi(t)^{i-1}\Psi(\chi(t))^j.
$$
Since $\chi(t+n)\to 0$, as $n\to +\infty$ and by (\ref{3.16}),  
$$\frac{\Psi(\chi(t+n))}{\chi(t+n)}\to\, 0,\,
\frac{\chi(t+1+n)}{\chi(t+n)}\to \, \lambda_1,\quad 
\text{ as} \quad n\to +\infty.$$
From $\Upsilon(t)=\chi(t)+\Psi(\chi(t))$, we have 
\begin{align}
\frac{\Upsilon(t+n+1)}{\Upsilon(t+n)}=
\frac{\chi(t+n+1)+\Psi(\chi(t+n+1))}{\chi(t+n)+\Psi(\chi(t+n))}&=
\frac{\frac{\chi(t+n+1)}{\chi(t+n)}+
\frac{\Psi(\chi(t+n+1))}{\chi(t+n+1)}\cdot \frac{\chi(t+n+1)}{\chi(t+n)}   }
{1+\frac{\Psi(\chi(t+n))}{\chi(t+n)}}\notag\\
&\to\, \lambda_1,\,\,
\text{as} \,\, n\,\to\, +\infty.\notag
\end{align}

 Conversely, if we put $\Upsilon(t)$ as (\ref{3.6}), 
where $\pi$ is an arbitrary periodic entire function, and 
$\Psi$ is a solution of (\ref{3.6}), 
then $\Upsilon(t)$ is a solution of (\ref{1.1}) such that 
$\Upsilon(t+n)\to 0$ as $n\to \,+\infty$. 
Furthermore then we have a solution $\chi$ of (\ref{3.3}) such that 
$$\Upsilon(t)=\chi(t)+\Psi(\chi(t)),$$
where $\chi(t+n)\to 0$ as $n\to +\infty$.
Hence we have 
$\frac{\Upsilon(t+1+n)}{\Upsilon(t+n)}\to \lambda_1$ as $n\to +\infty$. \par
\vspace{0.5cm}
Similarly in the proof of the above case, we can prove the case $m=2$ 
making use of the equations in (\ref{3.4}) and (\ref{3.5}). 
$\square$
\par
\vspace{0.7cm}



\appendix
\section*{Appendix A}
We put 
\begin{align}
p(t)&=\phi(p(t+1)+P_N(t+1),p(t+2)+P_N(t+2))-P_N(t)\notag\\
    &=g_1(t,p(t+1),p(t+2))+g_2(t)=g_3(t,p(t+1),p(t+2)),\notag
\end{align}

in which
\begin{align}
g_1(t,p(t+1),p(t+2))&=\phi(p(t+1)+P_N(t+1),p(t+2)+P_N(t+2) )\notag\\
  &\qquad-\phi(P_N(t+1),P_N(t+2))\notag \\ 
g_2(t)&=\phi(P_N(t+1),P_N(t+2))-P_N(t).\tag{1}
\end{align}
Since $\phi$ is holomorphic on $|w|\leqq \rho$, $|z|\leqq \rho$, using 
Cauchy's integral formula \cite{Ahlfors}, we have
$$
\frac{\partial \phi}{\partial w}(w,z)
=\frac{1}{2\pi i}\int_{|\xi|=\rho } 
\frac{\phi(\xi,z)}{(\xi-w)^2}d\xi.
$$
Therefore when $|w|\leqq \frac{\rho}{2}$, we have 
$|\xi-w|\geqq |\xi|-|w|\geqq \rho-\frac{\rho}{2}=\frac{\rho}{2}$ and 
$$
\Biggr |\frac{\partial \phi}{\partial w}(w,z)\Biggr |
\leqq \frac{1}{\pi}\int_{|\xi|=\rho}\frac{|\phi(\xi,z)|}{(\frac{\rho}{2})^2}|d\xi|
\leqq \frac{1}{\pi}\int_{|\xi|=\rho}\frac{K}{(\frac{\rho}{2})^2}|d\xi|=\frac{8K}{\rho}.
$$
When $|z|\leqq \frac{\rho}{2}$, similarly for $z$ we obtain 
$$\Biggr |\frac{\partial \phi}{\partial z}(w,z)\Biggr |\leqq \frac{8K}{\rho}.$$
Hence we have 
\begin{equation}
\Biggr |\frac{\partial \phi}{\partial w}\Biggr|,\, \Biggr |\frac{\partial \phi}{\partial z}\Biggr|
\leqq \frac{8K}{\rho} \qquad \mbox{for} \quad |w|,|z|\leqq \frac{\rho}{2}.\notag
\end{equation}

  Next we take $A$, and take $\eta$ sufficiently small such that $A\eta^{N+1}<\frac{\rho}{4}$. 
Then for sufficiently large $t$, we have 
     
$$|p(t)|\leqq A|\lambda^t|^{N+1}\leq A\eta^{N+1}<\frac{\rho}{4}.$$
And we have 
\begin{align}
|p(t+1)|&\leqq A|\lambda^{t+1}|^{N+1}
=A|\lambda|^{N+1}|\lambda^t|^{N+1}<\frac{\rho}{4},\notag\\
|p(t+2)|&\leqq A|\lambda^{t+2}|^{N+1}
=A|\lambda|^{2(N+1)}|\lambda^t|^{N+1}<\frac{\rho}{4}.\notag
\end{align}

Furthermore we can take $t$ so large that $|P_N(t+1)|,\, |P_N(t+2)|<\frac{\rho}{4}$, 
then we obtain 
$$|w|=|p(t+1)+P_N(t+1)|\leqq \frac{\rho}2,\,\,
|z|=|p(t+2)+P_N(t+2)|\leqq \frac{\rho}2.
$$
Since 
\begin{align}
g_1(t,p(t+1),p(t+2))
&=\int^1_0\frac{d}{dr}\phi(rp(t+1)+P_N(t+1),rp(t+2)+P_N(t+2))dr\notag\\
&=\int^1_0\{p(t+1)\frac{\partial \phi}{\partial w}(***)+
p(t+2)\frac{\partial \phi}{\partial z}(***)\}dr,\notag
\end{align}
where $(***)=(rp(t+1)+P_N(t+1),rp(t+2)+P_N(t+2))$, we have

\begin{align}
&|g_1(t,p(t+1),p(t+2))|\notag\\
&\leqq \int^1_0\{|p(t+1)|\Biggr| \frac{\partial \phi}{\partial w}(***)\Biggr |+
|p(t+2)|\Biggr| \frac{\partial \phi}{\partial z}(***)\Biggr |\}dr,\notag\\
&\leqq \int_0^1 \{A|\lambda^t|^{N+1}|\lambda|^{N+1}\cdot \frac{8K}{\rho}
+A|\lambda^t|^{N+1}|\lambda|^{2(N+1)}\cdot \frac{8K}{\rho}\}dr
\leqq \frac{16K}{\rho}A|\lambda|^{N+1}\cdot |\lambda^t|^{N+1}.\tag{2}
\end{align}
From definition of $P_N$ and (1), 
we have 
\begin{equation}
|g_2(t)|\leqq K_2|\lambda^t|^{N+1},\tag{3}
\end{equation}
with a constant $K_2$ which depends on $N$.
From (2) and (3), we have 
$$|T_1[p](t)|\leqq |g_1(t,p(t+1),p(t+2))|+|g_2(t)|
\leqq \Biggr (\frac{16K}{\rho}A|\lambda|^{N+1}+K_2\Biggr )|\lambda^t|^{N+1}.
$$
If we suppose $N$ is so large that 
$\frac{16K}{\rho}|\lambda|^{N+1}<\frac14,$ 
then we have 
$$|T_1[p](t)|
\leqq \Biggr(\frac{1}{4}A+K_2\Biggr)|\lambda^t|^{N+1}.
$$
Furthermore we take $A$ so large that 
$A>\frac 43 K_2,$ 
then  
\begin{equation}
|T_1[p](t)|<A|\lambda^t|^{N+1}.\notag
\end{equation}
So we obtain that $T_1$ in (\ref{2.8}) maps $J(A,\eta)$ into itself.\par

\section*{Appendix B}
Suppose there is another fixed point 
$p^*(t)=p_N^*(t)\in J(A^*,\eta^*)$. Put
\begin{gather}
A_0=\max(A,A^*), \quad \eta_0\leqq\min (\eta,\eta^*),\notag\\
u(t)=p_N(t)+P_N(t),\quad u^*(t)=p_N^*(t)+P_N(t),\notag
\end{gather}
and
$$q(t)=p_N^*(t)-p_N(t).$$
Then we have $|q(t)|\leqq 2A_0|\lambda^t|^{N+1}$. 
From (\ref{2.8}), we have
\begin{align}
q(t)&=\{\phi(p_N^*(t+1)+P_N(t+1),p_N^*(t+2)+P_N(t+2))-P_N(t)\}\notag\\
     &\qquad\qquad-\{\phi(p_N(t+1)+P_N(t+1),p_N(t+2)+P_N(t+2))-P_N(t)\}
\notag\\
&=\phi(q(t+1)+u_N(t+1),q(t+2)+u_N(t+2))-\phi(u_N(t+1),u_N(t+2))\notag\\
&=\int^1_0\{
q(t+1)\frac{\partial \phi}{\partial w}(****)
  +q(t+2)\frac{\partial \phi}{\partial z}(****)
\}dr\notag
\end{align}
where $(****)=(rq(t+1)+u_N(t+1),rq(t+2)+u_N(t+2))$. 
If $\eta_0$ is sufficiently small, then we have 
$$\Biggr|\frac{\partial \phi}{\partial w}(****)\Biggr|,\,\,
\Biggr| \frac{\partial \phi}{\partial z}(****) \Biggr|<\frac{8K_1}{\eta},$$
and we suppose $N$ is sufficiently large such that 
$|\lambda|^{N+1}<\frac {\rho}{64K_1}.$ 
Thus we have
\begin{align}
|q(t)|&\leqq \int^1_0 \frac{8K_1}{\rho}(|q(t+1)|+|q(t+2)|)dr\notag\\
&\leqq \int^1_0 \frac{8K_1}{\rho}|\lambda|^{N+1}
(2A_0|\lambda^t|^{N+1}+2A_0|\lambda^t|^{N+1})dr
<\frac12 A_0|\lambda^t|^{N+1}.\notag
\end{align}
Then
$$
|q(t)|=|p_N^*(t)-p_N(t)|\leqq \frac 12 A_0|\lambda^t|^{N+1}
=\Biggr(\frac14\Biggr)\cdot 2A_0|\lambda^t|^{N+1},\quad \mbox{for}\,\, t\in S(\eta_0).
$$

Next we consider $q(t)$ in which $|q(t)|\leqq \frac14\cdot 2A_0|\lambda^t|^{N+1}$ and 
repeat this procedure, then we have $|q(t)|\leqq (\frac14)^2 \cdot 2A_0|\lambda^t|^{N+1}$. 
  Repeating this procedure $k$ times we obtain
$$|p_N^*(t)-p_N(t)|<\Biggr( \frac 14\Biggr) ^k(2A_0)|\lambda^t|^{N+1},\quad
k=1,2,\cdots.$$
Letting $k\to \infty$, we have
$$p_N^*(t)=p_N(t),\quad t\in S(\eta_0).$$
\quad
 Thus $p_N^*(t)=p^*(t)$ and $p_N(t)=p(t)$ are holomorphic in
$|\lambda^t|\leqq \min(\eta,\eta^*)$ and $p^*(t)\equiv p(t)$ 
in $t\in S(\eta_0)$. Hence $p_N^*(t)=p_N(t)$ can be continued 
analytically to $S(\eta_1)$, $\eta_1=\max(\eta, \eta^*)$. 
$\square$\par

\section*{Appendix C}

Here we will show that the solution $u_N(t)$, given by $u_N(t)=p_N(t)+P_N(t)$ 
does not depend on $N$ in both cases in i).
  Let $p_N(t)\in J(A_N,\eta_N)$ and $p_{N+1}(t)\in J(A_{N+1},\eta_{N+1})$
 be fixed points of $T_1$, and
$$u_{N+1}(t)=p_{N+1}(t)+P_{N+1}(t)
          =p_{N+1}(t)+a_{N+1}\lambda^{(N+1)t}+P_N(t)
          =\tilde{p}_N(t)+P_N(t).
$$
\begin{align}
|\tilde{p}_N(t)|
=|p_{N+1}(t)+a_{N+1}\lambda^{(N+1)t}|
&\leqq A_{N+1}|\lambda^t|^{N+2}
+|a_{N+1}|\cdot |\lambda^t|^{N+1}\notag\\
&=(A_{N+1}|\lambda^t|+|a_{N+1}|)|\lambda^t|^{N+1}
=A_N^*|\lambda^t|^{N+1}\notag,
\end{align}
where $A^{*}=A_{N+1}|\lambda^t|+|a_{N+1}|$. 
We put $A=\max( A_N, A_N^{*})$. 
by uniqueness of fixed point, $\tilde{p}_N(t)=p_N(t)$ for
$t \in S(\eta_N)\cap S(\eta_{N+1})$. Thus
$$u_{N+1}(t)=u_{N}(t)\quad \mbox{in} \quad S(\eta_N)\cap S(\eta_{N+1}).$$
\quad By analytic prolongation \cite{Ahlfors}, both of $u_N(t)$ and $u_{N+1}(t)$ are holomorphic in
$S(\eta_N)\cap S(\eta_{N+1})$ and coincide there. Hence both of them are 
continued analytically  to $S(\eta_N)\cup S(\eta_{N+1})$ and 
$$u_{N+1}(t)=u_{N}(t)\quad \mbox{in} \quad S(\eta_N)\cup S(\eta_{N+1}). \,\,
\square $$ 

\end{document}